\documentclass[]{article}
\usepackage{graphicx,amsfonts,amsmath,amsthm,amssymb,amscd,longtable}

\def\tr{{\rm tr}}
\def\PSL{{\rm PSL}(2,{\Bbb C})}
\def\T{{\cal T}}

\newtheorem{theorem}{Theorem}[section]
\newtheorem{lemma}[theorem]{Lemma}

\newtheorem{cor}[theorem]{Corollary}
\newtheorem{remark}[theorem]{Remark}

\title{Discrete ${\cal RP}$ groups with a parabolic generator}
\author
{Elena Klimenko
\and
Natalia Kopteva
\thanks{The authors were partially supported by Gettysburg College
Research and Professional Development Grant, 2003--2004.}}

\date{\today}

\begin{document}

\maketitle

\begin{abstract}
We deal with two-generator subgroups of ${\rm PSL}(2,{\Bbb C})$
with real traces of both generators and their commutator.
We give discreteness criteria for these groups when at least 
one of the
generators is parabolic. We also present a list of the corresponding
orbifolds.
\end{abstract}


\noindent
Secondary: 20H10, 22E40, 57M60, 57S30.

\normalsize

\section{Introduction}

A two-generator subgroup $\Gamma=\langle f,g\rangle$ of $\PSL$
is determined up to conjugacy
by its parameters
$\beta(f)={\rm tr}^2f-4$, $\beta(g)={\rm tr}^2g-4$, and
$\gamma(f,g)={\rm tr}[f,g]-2$ whenever $\gamma(f,g)\not=0$
\cite{GM94-2}.

We are concerned with the class of $\cal RP$ {\it groups}
(two-generator groups with real parameters):
$$
{\cal RP}=\lbrace\Gamma|\Gamma=\langle f,g\rangle
{\rm \ for\ some\ }
f,g\in{\rm PSL}(2,{\Bbb C}) {\rm \ with\ } 
\beta(f),\beta(g),\gamma(f,g)\in {\Bbb R}\rbrace.
$$
Since discreteness conditions for elementary, Fuchsian, and NEC groups
are available, we consider only the non-elementary
${\cal RP}$ groups $\Gamma=\langle f,g\rangle$
without invariant planes.
In this paper, we deal with the case of $f$ elliptic, parabolic,
or hyperbolic (see the definitions in the beginning of Section~2) and
$g$ parabolic.
We give criteria for discreteness of these groups (Theorem~\ref{one_par}) 
and
for each discrete $\Gamma$ we obtain a presentation and the corresponding
Kleinian orbifold $Q(\Gamma)$ (Theorem~\ref{cor2e}).
Theorem~\ref{one_par} appeared in the PhD thesis of the first
author~\cite{Kli89-1}, but its proof was not published yet.

The groups generated by two parabolic elements were under study earlier.
In an unpublished preprint~\cite{Par89}, Parker gave a discreteness
criterion for an arbitrary group generated by two parabolic elements.
However, to apply this criterion a simple
condition needs to be checked 
for {\it all} elements of the group. In fact, the author
uses this criterion to give a nice discreteness condition for
Fuchsian groups (Proposition~3.1) and ${\cal RP}$ groups
(Proposition~3.2). The latter proposition and our Corollary~\ref{cor1}
for the case $\beta(f)=0$ imply each other.

Adams~\cite{Ada95} gives some necessary condition for discreteness
of a two parabolic generator non-elementary subgroup of $\PSL$.
He proves existence of universal upper bounds on the ``length''
of each of the parabolic generators and on the ``distance'' between them
(the lengths and distances are measured in a canonical choice of cusp
boundaries).
Using Thurston's Orbifold Theorem, he shows also that 
a non-elementary orientable finite volume hyperbolic 3-manifold
$M$ has fundamental group generated by two parabolic if and only if 
$M$ is the compliment
of a two-bridge link in ${\Bbb S}^3$ that is not a 2-braid.
Agol~\cite{Ago02} generalized Adams' argument to classify all
two parabolic generator {\it 3-orbifolds}.

\section{Discreteness criteria}

Recall that an element $f\in\PSL$ with real $\beta(f)$
is
{\it elliptic}, {\it parabolic}, {\it hyperbolic}, or {\it $\pi$-loxodromic}
according to whether
$\beta(f)\in[-4,0)$, $\beta(f)=0$, $\beta(f)\in(0,+\infty)$, or
$\beta(f)\in(-\infty,-4)$.
If $\beta(f)\notin[-4,+\infty)$,
then $f$ is called {\it strictly loxodromic}.
Among all
strictly loxodromic elements, only $\pi$-loxodromics have
real $\beta(f)$.

Every element $f\in{\rm PSL}(2,{\Bbb C})$ with real $\beta(f)$ has
invariant planes. The following lemma characterizes the non-elementary
${\cal RP}$ groups in terms of the invariant planes of generators.

\begin{lemma}\label{inv_planes}
Let $\Gamma=\langle f,g\rangle$ be a non-elementary subgroup of $\PSL$
with real $\beta(f)$ and $\beta(g)$. Then $\gamma(f,g)$ is real
if and only if either
\begin{itemize}
\item[$(1)$] $f$ and $g$ have a common invariant plane or
\item[$(2)$] each of the generators $f$ and $g$ has an invariant
plane orthogonal to all invariant planes of the other generator.
\end{itemize}
\end{lemma}

\noindent
{\it Proof} follows from Theorems~1--3 of \cite{KK02}.\qed

\begin{remark}
{\rm
Clearly, all non-elementary ${\cal RP}$ groups without invariant plane
satisfy the condition~(2) of Lemma~\ref{inv_planes}. 
Theorem~\ref{one_par} below characterizes all such discrete groups in case
when one of the generators is parabolic and the other has real trace.
}
\end{remark}

An elliptic element $f$ of order $n$ is said to be {\it non-primitive} 
if $f$ is a rotation through $2\pi k/n$, where $k$ and $n$ are coprime
($1<k<n/2$). If $f$ is a rotation through $2\pi/n$, then $f$ is called
{\it primitive}.

It is easy to see that if $f$ is a non-primitive elliptic
element of order $n$, then there exists an integer
$r\geq 2$ such that $f^r$
is a primitive elliptic element of the same order $n$.
It is clear that
$\langle f,g\rangle=\langle f^r,g\rangle$.
Therefore,
we assume without loss of generality that the elliptic
generator is primitive (see also Remark~\ref{rem26}).

\begin{theorem}\label{one_par}
Let $f\in \PSL$ be a hyperbolic, parabolic, or 
primitive elliptic element of order
$n\geq 3$,
let $g\in \PSL$ be a parabolic element,
and let $\Gamma=\langle f,g\rangle$ be a
non-elementary ${\cal RP}$ group without invariant plane.
Then:
\begin{itemize}
\item[$(1)$] there exists an element $h\in \PSL$ such that $h^2=fgf^{-1}g^{-1}$
and $(hg)^2=1$;
\item[$(2)$] $\Gamma=\langle f,g\rangle$ is discrete if and only if \
$h$ is a hyperbolic, parabolic, or primitive elliptic element of order
$p\geq 3$.
\end{itemize}
\end{theorem}

\noindent{\it Proof. }
We start with construction of a reflection group $\Gamma^*$
containing $\Gamma$
as a subgroup of finite index. Such a group is discrete if and
only if so is $\Gamma$.
Then we find discreteness criteria for $\Gamma^*$ and rewrite them
as simple conditions on the generators of~$\Gamma$.

\medskip
\noindent
{\bf 1.} {\it Construction of $\Gamma^*$ and a polyhedron $\cal T$
bounded by the planes of reflections of $\Gamma^*$.}
Let $f$ and $g$ be as in the statement of the theorem.
Since $\Gamma=\langle f,g\rangle$ is non-elementary, $f$ and $g$
have no fixed point in common.
Let $Q\in\partial {\Bbb H}^3$ be the fixed point of $g$
and let
$\zeta$ be the invariant plane of $f$ that passes through~$Q$.

Since $\Gamma=\langle f,g\rangle$ is a non-elementary
${\cal RP}$ group without
invariant plane,
it follows from Lemma~\ref{inv_planes} that
there is an invariant plane $\eta$ of $g$ 
orthogonal to
all invariant planes of $f$ and, in particular, to $\zeta$.
Note that if $f$ is elliptic, the axis of $f$ lies in $\eta$; and
if $f$ is hyperbolic, the axis of $f$ is orthogonal to $\eta$.

We see that there exists planes $\sigma$ and $\tau$ so that 
$f=R_\sigma R_\eta$
and $g=R_\tau R_\zeta$
(we denote the reflection in a plane $\kappa$ by $R_\kappa$).
Clearly, $\tau$ and $\zeta$ are parallel and meet at 
$Q\in\partial {\mathbb H}^3$.
Since $\eta$ is an invariant plane of $g$, $\eta$ is orthogonal to
both $\tau$ and $\zeta$.
Similarly, since $\zeta$ is an invariant plane of $f$,
$\zeta$ is orthogonal to both $\sigma$ and~$\eta$.

\begin{figure}[htbp]
\centering
\begin{tabular}{c}
\includegraphics[width=8 cm]{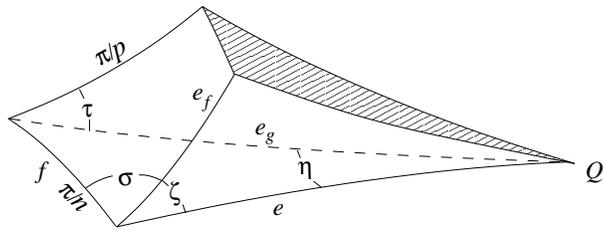}\\
\rule[-5ex]{0ex}{5ex}(a) $f$ is elliptic\\
\includegraphics[width=7 cm]{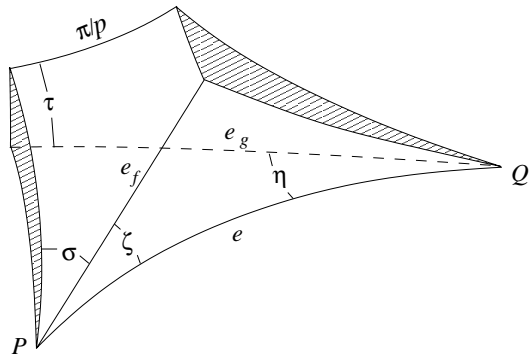}\\
\rule[-5ex]{0ex}{5ex}(b) $f$ is parabolic\\
\includegraphics[width=7 cm]{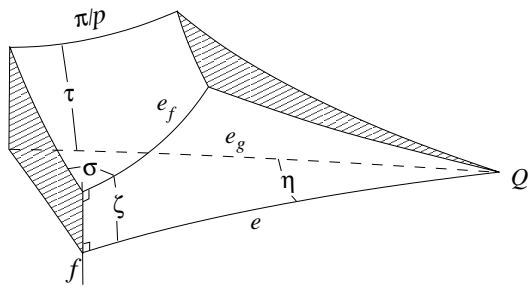}\\
(c) $f$ is hyperbolic
\end{tabular}
\caption{Fundamental polyhedra for $\Gamma^*$}\label{fund_poly}
\end{figure}

If $f$ is elliptic then the planes $\eta$ and $\sigma$
intersect at an angle of $\pi/n$; the line of their intersection
is the axis of~$f$ (see Figure~\ref{fund_poly}(a)).
If $f$ is parabolic then $\eta$ and $\sigma$ are parallel and
meet at a point $P\in\partial{\Bbb H}^3$ (see Figure~\ref{fund_poly}(b));
moreover, $\zeta$ passes through $P$.
If $f$ is hyperbolic, the planes $\eta$ and $\sigma$ are disjoint,
the axis of
$f$ is orthogonal to both $\eta$
and $\sigma$ (see Figure~\ref{fund_poly}(c)).

Consider half-turns $e=R_\eta R_\zeta=R_\zeta R_\eta$,
$e_f=R_\sigma R_\zeta$, and $e_g=R_\tau R_\eta$.
Then we have
\begin{eqnarray*}
f&=&R_\sigma R_\eta=(R_\sigma R_\zeta)(R_\zeta R_\eta)=e_fe\quad{\rm and}\\
g&=&R_\tau R_\zeta=(R_\tau R_\eta)(R_\eta R_\zeta)=e_ge.
\end{eqnarray*}

We define two finite
index extensions of the group $\Gamma=\langle f,g \rangle$  as follows:
$\widetilde{\Gamma}=\langle f,g,e\rangle$ and
$\Gamma^*=\langle f,g,e,R_\eta\rangle$.

It is easy to see that $\widetilde\Gamma=\Gamma\cup\Gamma e$.
If $e\in\Gamma$ then $\widetilde\Gamma=\Gamma$ and if
$e\notin\Gamma$ then $\Gamma$ is a subgroup of index~2 in
$\widetilde\Gamma$. As we will see, both possibilities are
realized. Since, moreover, $\widetilde\Gamma$ is the orientation
preserving subgroup of index~2 in $\Gamma^*$,
the groups $\Gamma$, $\widetilde\Gamma$,
and $\Gamma^*$ are either all discrete or all non-discrete.

It is clear that
$\Gamma^*=\langle R_\eta,R_\zeta,R_\sigma,R_\tau\rangle$.

Consider the infinite volume polyhedron $\T$ bounded by the planes
$\eta$, $\zeta$, $\sigma$, and $\tau$. As we have seen before,
$\eta$ and $\zeta$, $\eta$ and $\tau$, and $\zeta$ and $\sigma$
are orthogonal; $\zeta$ and $\tau$ are parallel; $\eta$ and $\sigma$
either intersect at an angle of $\pi/n$, or are parallel or disjoint
depending on the type of $f$.
The planes $\sigma$ and $\tau$ may either intersect, or be parallel
or disjoint (this depends on the type of $fgf^{-1}g^{-1}$ as we
will show in part 2 of the proof).

If $\sigma$ and $\tau$ intersect, then the dihedral angle
of $\T$ between these planes
(denote it by $\pi/p$ with $p$ not necessarily an integer)
is acute.
Indeed, there exists a hyperbolic plane $\kappa_1$ orthogonal
to $\zeta$, $\sigma$, and $\tau$; such a plane passes through~$Q$.
The planes $\zeta$, $\sigma$, and $\tau$ cut off a hyperbolic
triangle~$\Delta$ with angles $0$, $\pi/2$, and $\pi/p$
from $\kappa_1$
($\Delta$ is shaded in Figures~\ref{fund_poly}(a)--\ref{fund_poly}(c)).
Therefore, $\pi/p<\pi/2$.
We keep the notation $\pi/p$ taking $p=\infty$ or~$\overline\infty$
for the cases of parallel and disjoint $\sigma$ and $\tau$, respectively.
(We regard $\overline\infty>\infty>x$,
$x/\infty=x/\overline\infty=0$, $\infty/x=\infty$, and
$\overline\infty/x=\overline\infty$ for every positive real $x$.)

Similarly, if $\eta$, $\sigma$, and $\tau$ do not have a common 
point in ${\Bbb H}^3\cup\partial{\Bbb H}^3$, there exists a hyperbolic
plane $\kappa_2$ orthogonal to them. The planes $\kappa_1$
and $\kappa_2$ cut off a finite volume polyhedron $\overline\T$
from $\T$. By the Andreev theorem,
$\overline \T$  exists
in hyperbolic space for all $p>2$ \cite{And70, Vin85}; 
therefore, so does $\T$.
In fact, it is $\overline\T$ what is drawn in Figure~\ref{fund_poly};
moreover, it is drawn under assumption that
$p<\infty$ and, for Figure~\ref{fund_poly}(a), $1/p+1/n>1/2$.

\medskip\noindent
{\bf 2.} {\it Existence of $h$ and a sufficient discreteness condition
for $\Gamma$.}
It is clear that if
\begin{equation}\label{condition-i}
\text{$p$ is an integer ($p>2$), $\infty$, or $\overline\infty$,}
\end{equation}
then $\T$ and reflections $R_\eta$, $R_\zeta$, $R_\sigma$, and $R_\tau$
satisfy the hypotheses of the Poincar\'e theorem~\cite{EP94},
$\Gamma^*$ is discrete,
and $\T$ is its fundamental polyhedron.

Now we rewrite the condition~(\ref{condition-i}) via conditions on some
elements of $\PSL$.

Let us prove that there is a unique $h$ in $\PSL$ that satisfies both
$h^2=[f,g]$ and $(hg)^2=1$; moreover, $h=R_\sigma R_\tau$.

Since $\langle f,g\rangle$ is a non-elementary ${\cal RP}$ group,
$[f,g]$ is parabolic, or hyperbolic, or elliptic.
If $[f,g]$ is parabolic then it has only one square root; if $[f,g]$
is hyperbolic or elliptic, then $[f,g]$
has exactly two square roots $h$ and $\overline h$ in $\PSL$.
Namely, if $[f,g]$ is
hyperbolic, then one of the roots is hyperbolic and the other
is $\pi$-loxodromic. If $[f,g]$ is elliptic, then 
$h$ and $\overline h$ are both elliptic.

Let us show that if we take $h=R_\sigma R_\tau$ then $h^2=[f,g]$
and $(hg)^2=1$ hold. Indeed,
\begin{eqnarray*}
h^2&=&(R_\sigma R_\tau)^2=(R_\sigma R_\zeta R_\zeta R_\tau)^2=(e_fg^{-1})^2\\
&=&(e_fee_g)^2=(e_fe)(e_ge)(ee_f)(ee_g)=fgf^{-1}g^{-1}.
\end{eqnarray*}
Moreover, $hg=(R_\sigma R_\tau)(R_\tau R_\zeta)=R_\sigma R_\zeta=e_f$.
So, $(hg)^2=e_f^2=1$.

Let us explain now what $\overline h$ is. If $[f,g]$ is hyperbolic,
then $\overline h$ is $\pi$-loxodromic with the same axis and translation
length as $h$. If $[f,g]$ is elliptic, then $\overline h$ is elliptic
with the same axis as $h$ and with rotation angle $(\pi-2\pi/p)$, while
$h$ is a rotation through $2\pi/p$ in the opposite direction.
It is clear that in both cases $(\overline hg)^2\not=1$.

This means that the only element $h$ that satisfies
both $h^2=[f,g]$ and $(hg)^2=1$ can be written as $h=R_\sigma R_\tau$.
Thus, part (1) of Theorem~\ref{one_par}
is proved.

The element
$h=R_\sigma R_\tau$ is a primitive elliptic element of order $p\geq 3$
if and only if the dihedral angle of $\cal T$ at the
edge $\sigma\cap\tau$ is equal to $\pi/p$, $p\in {\Bbb Z}$;
$h$ is parabolic (hyperbolic) if and only if $\sigma$ and
$\tau$ are parallel (disjoint, respectively).

Therefore, we have proved that
the condition~(\ref{condition-i}) is equivalent to the condition that
\begin{equation}\label{cond_h}
\text{$h$ is a hyperbolic, parabolic, or primitive elliptic
element of order $p>2$.}
\end{equation}
So, (\ref{cond_h}) implies that $\Gamma$ is discrete.

\medskip\noindent
{\bf 3.} {\it The sufficient condition (\ref{cond_h}) is also a
necessary condition.}
Now suppose that $\Gamma$ is discrete but (\ref{cond_h})
fails.
This means that $h$ is a non-primitive elliptic element of finite order,
i.e., $p=q/k>2$, where $q$ and $k$ are coprime, $k\geq 2$.
Let us show that this is impossible.

Consider the hyperbolic plane $\kappa_1$
that is orthogonal to the planes $\zeta$, $\sigma$,
and $\tau$. Since $\langle e_f,g\rangle\subset\widetilde\Gamma$
keeps $\kappa_1$ invariant and preserves orientation of $\kappa_1$,
 $\langle e_f,g\rangle$
acts on $\kappa_1$ as a subgroup of ${\rm PSL}(2,{\Bbb R})$.
However, $\langle e_f,g\rangle$ is not discrete
if $h=e_fg^{-1}$ is non-primitive elliptic, which is a rotation
through $2k\pi/q$ in our case, by \cite{Kna68} or \cite{Mat82}.

Theorem~\ref{one_par} is proved.
\qed

\begin{remark}
{\rm
For the benefit of the reader, we give a description of $\T$ in the
upper half-space model of hyperbolic 3-space
${\Bbb H}^3=\{(z,t)\,:\,z\in{\Bbb C}, t>0\}$
with the Poincar\'e metric
$ds^2=(|dz|^2+dt^2)/t^2$.

It suffices to assume that in
the proof of Theorem~\ref{one_par}, $\Gamma=\langle f,g \rangle$
is normalized so that
$
g=\left(
\begin{array}{cc}
1 & 1\\
0 & 1
\end{array}
\right)
$
and the fixed points of $f$ are $z_0$ and $-z_0$, $z_0\in{\Bbb C}$.
(If $f$ is parabolic, then the only fixed point of $f$ is $z_0=0$.)

Suppose that $f$ is elliptic. Since $\Gamma$ is a non-elementary
${\cal RP}$ group without invariant plane, both fixed points of $f$
lie in an invariant plane of $g$ by Lemma~\ref{inv_planes}.
Taking into account the fact that every invariant plane of $g$
is given by $\{(z,t)\in{\Bbb H}^3:{\rm Im}z=const\}$, we conclude that
$z_0=x_0$ is real.

Analogously, if $f$ is hyperbolic, the fixed points of $f$ are symmetric
to each other with respect to an invariant plane of $g$ and hence
$z_0=iy_0$, $y_0\in{\Bbb R}$.

So we have normalized $\Gamma$ so that $g(\infty)=\infty$, $g(0)=1$,
$f(z_0)=z_0$, $f(-z_0)=-z_0$, where $z_0$ equals to $x_0$, $0$, or
$iy_0$ ($x_0$ and $y_0$ are some non-zero real numbers), if $f$
is elliptic, parabolic, or hyperbolic, respectively.

Then the planes $\eta$, $\zeta$, $\sigma$, and $\tau$ in the proof are
given by
\begin{eqnarray*}
\eta&=&\{(z,t)\in{\Bbb H}^3:{\rm Im}z=0\};\\
\zeta&=&\{(z,t)\in{\Bbb H}^3:{\rm Re}z=0\};\\
\sigma&=&\{(z,t)\in{\Bbb H}^3:|z-iy_\sigma|^2+t^2=r^2\}, \
y_\sigma\in{\Bbb R}\backslash\{0\},\ r^2=y_\sigma^2+z_0^2;\\
\tau&=&\{(z,t)\in{\Bbb H}^3:{\rm Re}z=1/2\}.
\end{eqnarray*}

Moreover, $\kappa_1=\{(z,t)\in{\Bbb H}^3\,:\,{\rm Im}z=y_\sigma\}$
is the plane that plays a key role in part~3 of the proof.
It is clear that this plane is orthogonal to all $\zeta$, $\sigma$,
and $\tau$.
}
\end{remark}

\begin{cor}\label{cor1}
Let $f,g\in\PSL$, $\beta(f)\in[0,+\infty)$ or $\beta(f)=-4\sin^2(\pi/n)$,
$n\in{\Bbb Z}$, $n\geq 3$,
and let
$\beta(g)=0$. Suppose that $\gamma(f,g)<0$.
Then $\Gamma=\langle f,g\rangle$
is discrete if and only if one of the following holds:
\begin{enumerate}
\item $\gamma(f,g)\in(-\infty;-4]$;
\item $\gamma(f,g)=-4\cos^2(\pi/p)$, $p\in{\Bbb Z}$, $p\geq 3$.
\end{enumerate}
\end{cor}

\noindent
{\it Proof.} Since $\langle f,g\rangle$ is an ${\cal RP}$ group,
$\beta(g)=0$, and $f$ is not $\pi$-loxodromic,
$\gamma$ is a non-elementary group without invariant plane
if and only if
$\gamma(f,g)<0$ \cite[Theorem~4]{KK02}. So it is clear that the
hypotheses of Corollary~\ref{cor1} are equivalent to
those of Theorem~\ref{one_par}.

Therefore, to prove the corollary it suffices to rewrite
part (2) of Theorem~\ref{one_par} in terms of $\gamma(f,g)$.
Since $\gamma(f,g)=\tr[f,g]-2$ and $[f,g]=h^2$,
it is not difficult to find $\gamma(f,g)$.

The element $h$ is hyperbolic if and only if the planes
$\sigma$ and $\tau$ are disjoint.
(We denote planes
 as in the proof of Theorem~\ref{one_par}.)
Let $d$ be the hyperbolic distance between them. Since
$[f,g]=h^2=(R_\sigma R_\tau)^2$,
$$
\gamma(f,g)=\tr[f,g]-2=-2\cosh(2d)-2<-4
$$
(we must take $\tr[f,g]$ to be negative, because $\gamma(f,g)$
is negative by assumption).

The element $h$ is parabolic if and only if $[f,g]$ is parabolic, that is,
$\tr[f,g]=-2$
($\tr[f,g]=2$ would give $\gamma(f,g)=0$ which is impossible in our case).
Hence, $\gamma(f,g)=\tr[f,g]-2=-4$.

Thus, $h$ is hyperbolic or parabolic if and only if
$\gamma(f,g)\in(-\infty,-4]$, and part~1 of Corollary~\ref{cor1}
is proved.

\medskip
Now suppose that $h$ is an elliptic element with rotation angle $\varphi$,
where $\varphi/2=\pi/p<\pi/2$ is the dihedral angle of $\cal T$ between
$\sigma$ and $\tau$. Then $[f,g]=h^2$ is also elliptic with rotation angle
$2\varphi$. Since $\tr[f,g]$ is well-defined (does not depend on the choice
of representatives for $f$ and $g$ in ${\rm SL}(2,{\Bbb C})$) we can
determine which formula, $\tr[f,g]=+2\cos\varphi$
or $\tr[f,g]=-2\cos\varphi$, is correct. The easiest way to do this
is by using the continuity of $\tr[f,g]$ as a function of $\varphi$ and
the limit condition $\tr[f,g]\to -2$ as $\varphi\to 0$.
So we must take $\tr[f,g]=-2\cos\varphi$ where $\varphi<\pi$ is the
doubled dihedral angle of $\cal T$.

Vice versa, if $\tr[f,g]$ is given, we can use the formula
$\tr[f,g]=-2\cos\varphi$, $\varphi<\pi$, 
to determine the rotation angle $\varphi$
of the element $h$ from Theorem~\ref{one_par}.

Thus, $h$ is a primitive elliptic element of order $p$, that is
$\varphi=2\pi/p$, if and only if
$$
\gamma(f,g)=\tr[f,g]-2=-2\cos(2\pi/p)-2=-4\cos^2(\pi/p), \ p\in{\Bbb Z}.
$$

Corollary~\ref{cor1} is proved.
\qed

\begin{remark}\label{rem26}
{\rm
For simplicity, in the statement of Corollary~\ref{cor1}
the elliptic generator $f$ is assumed to be primitive.
If $f$ is non-primitive elliptic then Corollary~\ref{cor1}
still can be used to verify
whether $\Gamma$ is discrete, but first we must replace
the triple $(\beta(f),\beta(g),\gamma(f,g))$,
where $\beta(f)=-4\sin^2(q\pi/n)$, $(q,n)=1$, $1<q<n/2$,
with the new triple $(\widetilde\beta,\beta(g),\widetilde\gamma)$,
where $\widetilde\beta=-4\sin^2(\pi/n)$ and
$\widetilde\gamma=(\widetilde\beta/\beta(f))\gamma(f,g)$.
The new triple corresponds to the same group by
Gehring and Martin~\cite{GM94-1} (cf. \cite[Remark~2, p.~262]{KK02}).
}
\end{remark}

\section{Kleinian orbifolds with a parabolic generator}

Let $\Gamma$ be a non-elementary Kleinian group.
Denote by $\Omega(\Gamma)$ the discontinuity set of $\Gamma$.
Following \cite{BP00}, we say that the {\it Kleinian orbifold}
$Q(\Gamma)=({\Bbb H}^3\cup\Omega(\Gamma))/\Gamma$
is an orientable $3$-orbifold with a complete hyperbolic structure
on its interior ${\Bbb H}^3/\Gamma$ and a conformal structure
on its boundary $\Omega(\Gamma)/\Gamma$.

We use the following notations:

\begin{itemize}
\item
$GT[n,m;q]=\langle f,g\,|\,f^n=g^m=[f,g]^q=1\rangle,$
where $n,m,q\in\{2,3,\dots\}\cup\{\infty,\overline\infty\}$.
If some relation has power $\overline\infty$, 
then we simply remove this relation from the presentation.
Further, if some relation has power $\infty$ and we keep 
it, we get a Kleinian group presentation. To obtain an abstract group
presentation, we need to remove this relation as well.

For example, the group $GT[n,\infty;\overline\infty]$
has Kleinian group presentation
$\langle f,g\,|\,f^n=g^\infty=1\rangle$ and is isomorphic to
${\Bbb Z}_n*{\Bbb Z}$.

\item
$Tet[n,m;q]=\langle x,y,z\,|\,
x^2=y^2=z^n=(xy^{-1})^m=(yz^{-1})^2=(zx^{-1})^q=1\rangle$,
where $n,m,q\in{\Bbb Z}\cup\{\infty,\overline\infty\}$.
For finite $n$, $m$, and $q$, this group
is a tetrahedron group generated by rotations in edges of a
face of an orthoscheme with the Coxeter diagram
$\circ\overset{n}{-}\circ\overset{q}{-}\circ\overset{m}{-}\circ$.
Note that
$Tet[n,m;q]\cong Tet[m,n;q]$.
\end{itemize}

\begin{theorem}\label{cor2e}
Let $\Gamma=\langle f,g\rangle$ be a non-elementary discrete 
${\cal RP}$ group
without invariant plane. Let $\beta(g)=0$ and
let $\beta(f)=-4\sin^2(\pi/n)$,
$n\in{\Bbb Z}$, $n\geq 3$, or $\beta(f)\in[0,+\infty)$. 
Put $n=\infty$ for
$\beta(f)=0$ and $n=\overline\infty$ for $\beta(f)\in(0,+\infty)$.
\begin{enumerate}
\item If $\gamma(f,g)\in(-\infty;-4)$, then $\Gamma$ is isomorphic to
$GT[n,\infty;\overline\infty]$.
\item If $\gamma(f,g)=-4$, then $\Gamma$ is isomorphic to
$GT[n,\infty;\infty]$.
\item If $\gamma(f,g)=-4\cos^2(\pi/p)$, $(p,2)=2$,
$p\geq 4$, then
$\Gamma$ is isomorphic to
$GT[n,\infty;p/2]$.
\item If $\gamma(f,g)=-4\cos^2(\pi/p)$, $(p,2)=1$,
$p\geq 3$, then $\Gamma$ is isomorphic to
$Tet[n,\infty;p]$.
\end{enumerate}
\end{theorem}

\begin{figure}[htbp]
\centering
\begin{tabular}{ccc}
\multicolumn{3}{c}%
{$\pi_1^{orb}(Q)\cong\langle f,g\,|\,f^n=g^\infty=[f,g]^q=1\rangle$, %
$n\geq 3$, $q\geq 2$}\\
\rule[-4ex]{0ex}{6ex}$q<\infty$ \quad & \quad $q=\infty$ \quad & \quad $q=\overline\infty$\\
\includegraphics[height=2.3 cm]{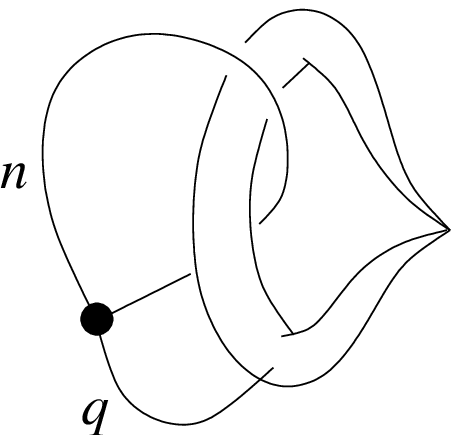}\quad & \quad
\includegraphics[height=2.3 cm]{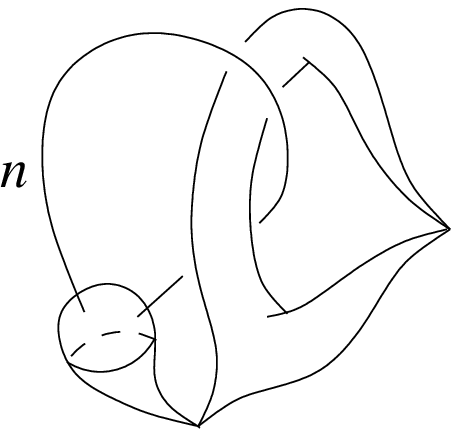}\quad & \quad
\includegraphics[height=2 cm]{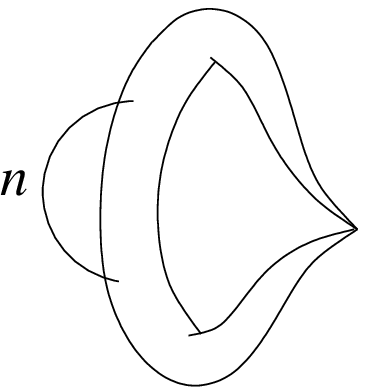}\\
\rule{0ex}{5ex}(1-a) & (1-b) & (1-c)\\

\multicolumn{3}{c}%
{\rule{0ex}{5ex}$\pi_1^{orb}(Q)\cong\langle f,g\,|\,f^\infty=g^\infty=[f,g]^q=1\rangle$, %
$q\geq 2$}\\
\rule[-4ex]{0ex}{6ex}$q<\infty$ \quad & \quad $q=\infty$ \quad & \quad $q=\overline\infty$\\
\includegraphics[height=2.3 cm]{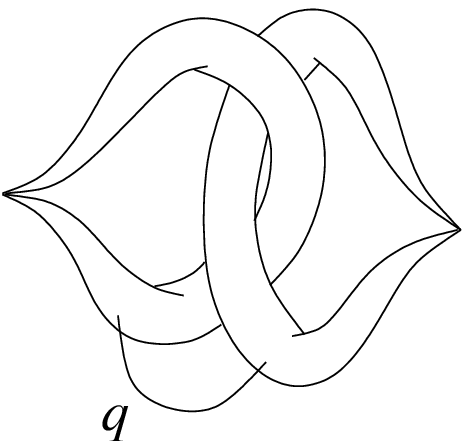}\quad & \quad
\includegraphics[height=2.3 cm]{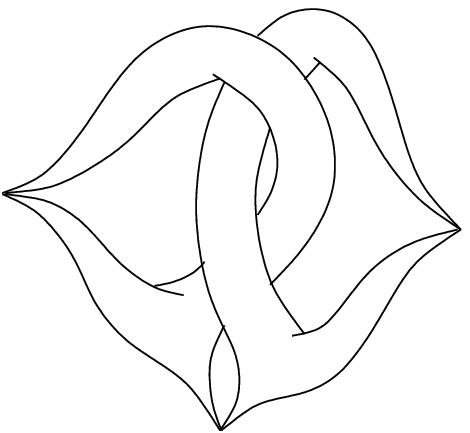}\quad & \quad
\includegraphics[height=1.5 cm]{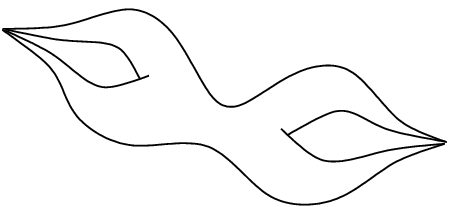}\\
\rule{0ex}{5ex}(2-a) & (2-b) & (2-c)\\

\multicolumn{3}{c}%
{\rule{0ex}{5ex}$\pi_1^{orb}(Q)\cong\langle f,g\,|\,g^\infty=[f,g]^q=1\rangle$, %
$q\geq 2$}\\
\rule[-4ex]{0ex}{6ex}$q<\infty$ \quad & \quad $q=\infty$ \quad & \quad $q=\overline\infty$\\
\includegraphics[height=2.3 cm]{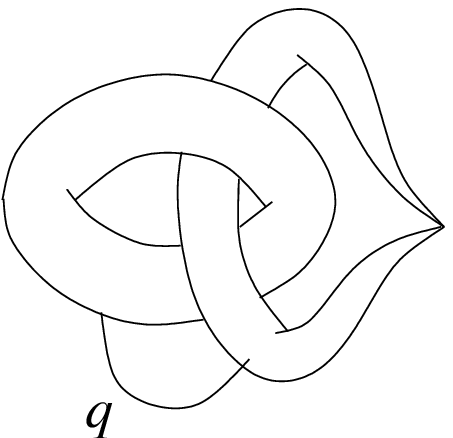}\quad & \quad
\includegraphics[height=2.3 cm]{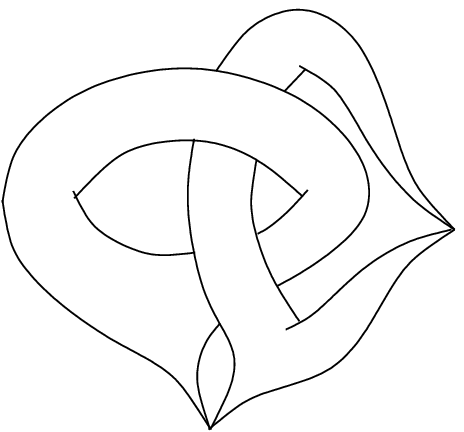}\quad & \quad
\includegraphics[height=1.5 cm]{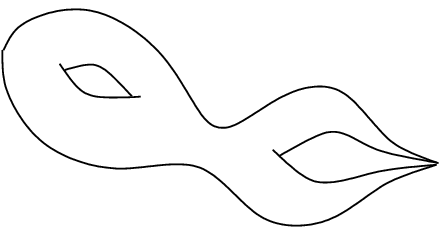}\\
\rule{0ex}{5ex}(3-a) & (3-b) & (3-c)
\end{tabular}
\caption{Orbifolds $Q$ with fundamental group $GT[n,\infty;q]$}\label{gt_orb}
\end{figure}

\begin{figure}[htbp]
\centering
\begin{tabular}{ccc}
\multicolumn{3}{c}%
{$\pi_1^{orb}(Q)\cong\langle f,g,e\,|\,f^n=g^\infty=e^2=(fe)^2=(ge)^2=
(gfe)^p=1\rangle$, $p\geq 3$}\\
\rule{0ex}{5ex}$f$ is elliptic \quad & \quad $f$ is parabolic \quad
& \quad $f$ is hyperbolic\\
\rule[-4ex]{0ex}{6ex}$n<\infty$\quad & \quad $n=\infty$ \quad
& \quad $n=\overline\infty$\\
\includegraphics[height=1.5 cm]{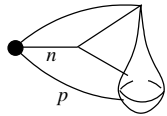} \quad & \quad
\includegraphics[height=1.5 cm]{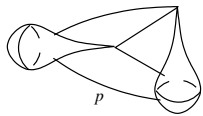} \quad & \quad
\includegraphics[height=1.5 cm]{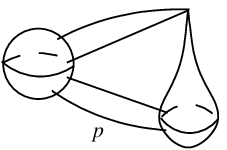}\\
\rule{0ex}{5ex}(a) \quad & \quad (b) \quad & \quad (c)
\end{tabular}
\caption{Orbifolds $Q$ with fundamental group $Tet[n,\infty;p]$}\label{tet_orb}
\end{figure}

\begin{remark}
{\rm
For all $\Gamma$ from Theorem~\ref{cor2e}, the orbifolds $Q(\Gamma)$
are shown in Figures~\ref{gt_orb} and~\ref{tet_orb} where the singular
sets and boundaries of $Q(\Gamma)$ are drawn. 
The indices on edges correspond to the orders of cone points, 
indices 2 are omitted.
Each $Q(\Gamma)$
is embedded in ${\Bbb S}^3={\Bbb R}^3\cup\{\infty\}$
so that $\infty$ is a non-singular interior point of $Q(\Gamma)$.

Since we obtain an orbifold by gluing faces of a fundamental polyhedron
for the group action on ${\Bbb H}^3$, not only topological, but also
metric structure is uncovered (the lengths of singular geodesics,
the structure of cusps, etc.).
In fact, since all fixed points of parabolic elements of $\Gamma$
belong to the limit set $\Lambda(\Gamma)$, they have no images
in $Q(\Gamma)$. For example, in Figure~\ref{gt_orb}(2-b),
the boundary of $Q(\Gamma)$ is the union of two thrice-punctured
2-spheres.

The fat vertices in Figures~\ref{gt_orb} and~\ref{tet_orb}
are either singuler points of the orbifold, or punctures, or correspond
to removed open balls.
In the first case the vertex
corresponds to the image of a point of ${\mathbb H}^3$,
in the second it corresponds to the image of a point of the
limit set $\Lambda(\Gamma)$ and, therefore, does not belong
to the orbifold, and in the last case the fat vertex corresponds 
to a boundary component of $\Omega(\Gamma)/\Gamma$.The type of a fat vertex
depends on the indices at the edges incident to it. For example, in
Figure~\ref{tet_orb}(a), the fat vertex is an interior singular point of
the orbifold if $1/2+1/n+1/p>1$, is a puncture if $1/2+1/n+1/p=1$,
and is a removed open ball if $1/2+1/n+1/p<1$.
}
\end{remark}

\noindent
{\it Proof of Theorem~\ref{cor2e}.}
We give a proof only for the case $\beta(f)=-4\sin^2(\pi/n)$,
$n\in{\Bbb Z}$. 

All parameters for the discrete groups in the statement of
Theorem~\ref{cor2e} are described in Corollary~\ref{cor1}. We will obtain
a presentation for each case by using the Poincar\'e polyhedron theorem.

\begin{figure}[htbp]
\centering
\begin{tabular}{c}
\includegraphics[width=7 cm]{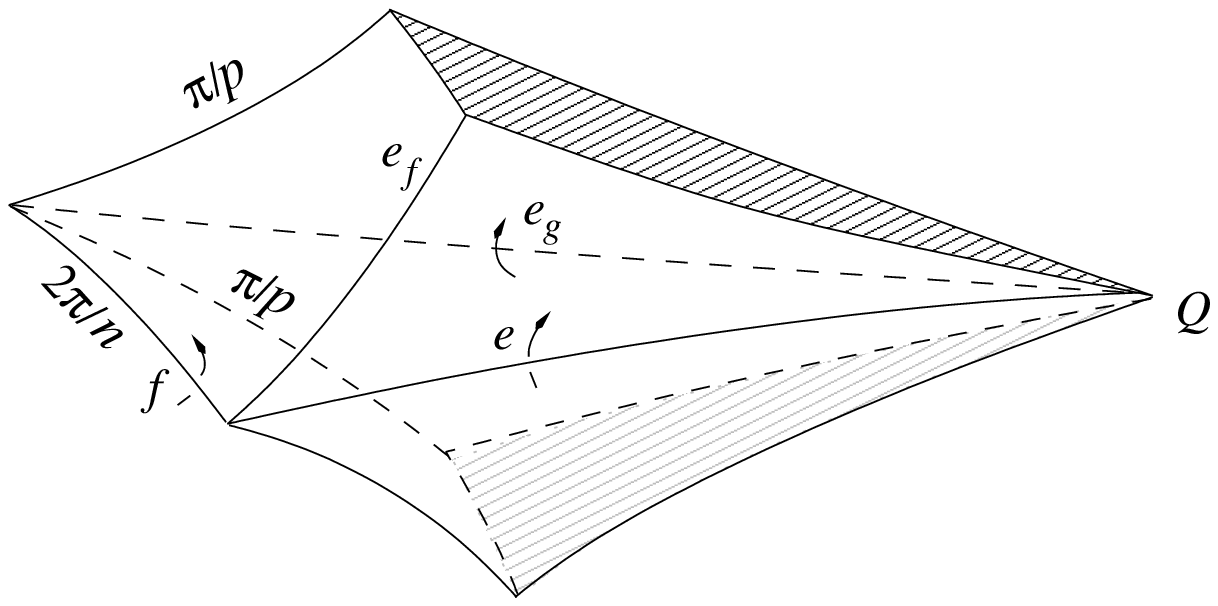}\\
(a)\\
\includegraphics[width=7 cm]{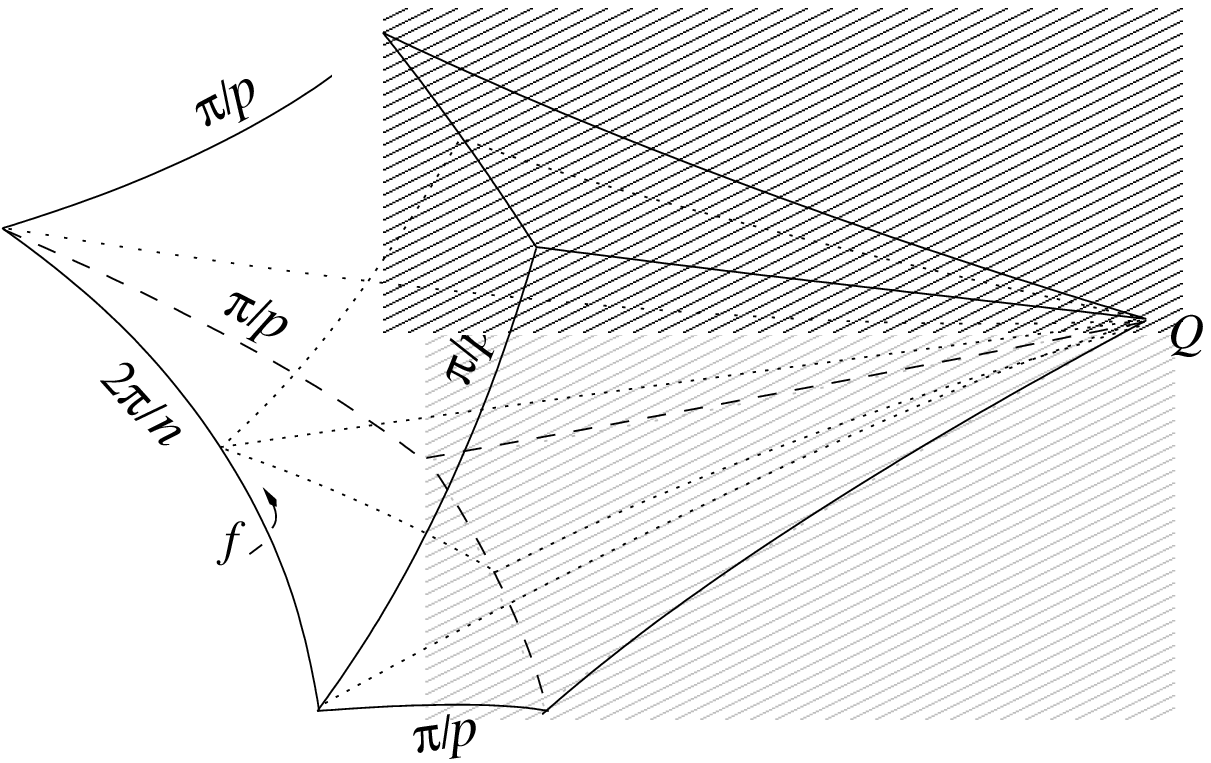}\\
(b)
\end{tabular}
\caption{Fundamental polyhedra for $\widetilde\Gamma$ and
$\Gamma$}\label{poincare_th}
\end{figure}

We start with construction of a fundamental polyhedron and a presentation for
the group $\widetilde\Gamma$ defined in the proof of Theorem~\ref{one_par}.
Since $\widetilde\Gamma$ is the orientation preserving index~2 subgroup
in $\Gamma^*$ and $\cal T$ is a fundamental polyhedron for $\Gamma^*$,
a fundamental polyhedron ${\cal P}$ for $\widetilde\Gamma$
consists of two copies
of~$\cal T$ (see Figure~\ref{poincare_th}(a)).

By applying the Poincar\'e polyhedron theorem to ${\cal P}$ and face pairing
transformations $e$, $e_g$, and $f$, we get
$$
\widetilde\Gamma=\langle e,e_g,f\,|\,e^2=e_g^2=f^n=(e_ge)^\infty=(fe)^2=
(fe_g)^p=1\rangle,
$$
where $p$ is an integer, $\infty$, or $\overline\infty$.
Since $g=e_ge$,
$$
\widetilde\Gamma=\langle f,g,e\,|\,f^n=g^\infty=e^2=(fe)^2=(ge)^2=
(fge)^p=1\rangle.
$$
Note that if $p\geq 3$ is odd, then using the relations
$(fe)^2=(ge)^2=1$, from $(fge)^p=1$ we obtain
$e=(fgf^{-1}g^{-1})^{(p-1)/2}fg$. Hence, in this case
$\widetilde\Gamma=\Gamma$ and $\Gamma\cong Tet[n,\infty;p]$.
Identifying the faces of ${\cal P}$, we get the orbifold
${\Bbb H}^3/\Gamma$ shown in Figure~\ref{tet_orb}(a).

If $p\geq 4$ is even, $\infty$, or $\overline\infty$,
then $\Gamma$ is a subgroup of index~2 in $\widetilde\Gamma$.
To see this we apply the Poincar\'e theorem to a polyhedron 
consisting of four copies of~$\cal T$
(see Figure~\ref{poincare_th}(b)). Then
$$
\Gamma=\langle f,g\,|\,f^n=g^\infty=(fgf^{-1}g^{-1})^{p/2}=1\rangle.
$$
The orbifold $Q(\Gamma)$ for finite $q=p/2$ is shown in
Figure~\ref{gt_orb}(1-a).
For the case of parallel $\sigma$ and $\tau$ ($p=\infty$),
$Q(\Gamma)$ is shown in Figure~\ref{gt_orb}(1-b),
and for the case of disjoint $\sigma$ and $\tau$,
$Q(\Gamma)$
is shown in Figure~\ref{gt_orb}(1-c).

The cases when $f$ is parabolic or hyperbolic have similar proofs and we
leave them to the reader.
\qed

\noindent
Gettysburg College\\
Mathematics Department\\
300 North Washington St., CB 402\\
Gettysburg, PA 17325\\
USA\\
{\tt yklimenk@gettysburg.edu}\\
\\
CMI, Universit\'e de Provence\\
39, rue F. Joliot Curie,\\
13453 Marseille cedex 13\\
France\\
{\tt kopteva@cmi.univ-mrs.fr}

\end{document}